\documentclass[11pt]{article}

\usepackage{amssymb,amsthm,amscd,array}
\usepackage{graphics}
\usepackage[final]{epsfig}	

\newcommand{\bbR}{\mathbb{R}}

\newcommand{\bbC}{\mathbb{C}}

\newcommand{\Vect}{\mathrm{Vect}}
\newcommand{\CVect}{\mathrm{CVect}}
\newcommand{\TVect}{\mathrm{TVect}}

\newcommand{\cB}{{\mathcal{B}}}
\newcommand{\cF}{{\mathcal{F}}}
\newcommand{\cI}{{\mathcal{I}}}
\newcommand{\cH}{{\mathcal{H}}}

\newcommand{\fh}{\mathfrak{h}}

\chardef\s=110
\chardef\g=103

\begin{document}

\newtheorem{thm}{Theorem}[section]
\newtheorem{lem}[thm]{Lemma}
\newtheorem{cor}[thm]{Corollary}
\newtheorem{prop}[thm]{Proposition}
\newtheorem{rmk}[thm]{Remark}
\newtheorem{exe}[thm]{Example}
\newtheorem{defi}[thm]{Definition}

\def\a{\alpha}
\def\b{\beta}
\def\d{\delta}
\def\g{\gamma}
\def\om{\omega}
\def\r{\rho}
\def\s{\sigma}
\def\vfi{\varphi}
\def\vr{\varrho}
\def\l{\lambda}
\def\m{\mu}

\title{Vector fields in the presence of a contact structure}

\author{V. Ovsienko
\thanks{
CNRS, 
Institut Camille Jordan
Universit\'e Claude Bernard Lyon 1,
21 Avenue Claude Bernard,
69622 Villeurbanne Cedex,
FRANCE;
ovsienko@igd.univ-lyon1.fr
}}

\date{}

\maketitle

{\abstract{
We consider the Lie algebra of all vector fields on a contact manifold as a 
module over the Lie subalgebra of contact vector fields. 
This module is split into a direct sum of two submodules: the contact algebra itself and
the space of tangent vector fields. We study the geometric nature of these two modules.
}}

\thispagestyle{empty}

\section{Introduction}

Let $M$ be a (real) smooth manifold and $\Vect(M)$ the Lie algebra of all smooth vector fields
on $M$ with complex coefficients. We consider the case when $M$ is $(2n+1)$-dimensional and
can be equipped with a contact structure. For instance, if $\dim{}M=3$, and $M$ is compact
and orientable, then the famous theorem of 3-dimensional topology states that there is always
a contact structure on $M$.

Let $\CVect(M)$ be the Lie algebra of smooth vector fields on $M$ preserving the contact
structure. This Lie algebra naturally acts on $\Vect(M)$ (by Lie bracket). We will study the
structure of $\Vect(M)$ as a $\CVect(M)$-module.
First, we observe that $\Vect(M)$ is split, as a $\CVect(M)$-module, into a direct sum of two
submodules:
$$
\Vect(M)\cong\CVect(M)\oplus\TVect(M)
$$
where $\TVect(M)$ is the space of vector fields tangent to the contact distribution.
Note that the latter space is a $\CVect(M)$-module but not a Lie subalgebra of $\Vect(M)$.

The main purpose of this paper is to study the two above spaces geometrically.
The most important notion for us is that of \textit{invariance}. 
All the maps and isomorphisms we consider are invariant with respect to the group of contact
diffeomorphisms of $M$.
Since we consider only local maps, this is
equivalent to the invariance with respect to the action of the Lie algebra $\CVect(M)$.

It is known, see \cite{Ovs,OT}, that the adjoint action of $\CVect(M)$ has the following
geometric interpretation:
$$
\CVect(M)\cong\cF_{-\frac{1}{n+1}}(M),
$$
where $\cF_{-\frac{1}{n+1}}(M)$ is the space of (complex valued) tensor densities of degree
$-\frac{1}{n+1}$ on $M$, that is, of sections of the line bundle
$$
\left(\wedge^{2n+1}T^*_\bbC{}M\right)^{-\frac{1}{n+1}}
\to{}M.
$$
In particular, this provides the
existence of a nonlinear invariant functional on $\CVect(M)$ defined on the contact vector
fields with nonvanishing contact Hamiltonians.

The analogous interpretation of $\TVect(M)$ is more complicated:
$$
\TVect(M)\cong\Omega^2_0(M)\otimes\cF_{-\frac{2}{n+1}}(M),
$$
where $\Omega^2_0(M)$ is the space of 2-forms on $M$ vanishing on the contact distribution.
Here and below the tensor products are defined over $C^\infty(M)$.

We study the relations between $\TVect(M)$ and $\CVect(M)$.
We prove the existence of a non-degenerate skew-symmetric invariant bilinear
map
$$
\cB:\TVect(M)\wedge\TVect(M)\to\CVect(M)
$$
that measures the non-integrability, i.e., the failure of the Lie bracket of two tangent
vector fields to remain tangent.

In order to provide explicit formul{\ae}, we introduce a notion of Heisenberg structure
on $M$. 
Usually, to write explicit formul{\ae} in contact geometry, one uses the Darboux
coordinates. 
However, this is not the best way to proceed
(as already noticed in \cite{Kir}).
The Heisenberg structure provides a universal expression for a contact vector field and its
actions. 

\section{Contact and tangent vector fields}

In this section we recall the basic definitions of contact geometry.
We then prove our first statement on a decomposition of the Lie algebra of all smooth vector
fields viewed as a module over the Lie algebra of contact vector fields.

\subsection{Main definitions}

Let $M$ be a $(2n+1)$-dimensional manifold. A contact structure on $M$ is a codimension
1 distribution $\xi$ which is completely non-integrable. The distribution $\xi$ can be
locally defined as the kernel of a differential 1-form $\a$ defined up to multiplication by a
nonvanishing function.  Assume that $M$ is orientable, then the form $\a$ can be globally defined on
$M$. Complete
non-integrability means that
\begin{equation}
\label{CompNonint}
\Omega:=\a\wedge(d\a)^{n}
\neq0
\end{equation}
everywhere on $M$. In other words, $\Omega$ is a volume form.
The above condition is also equivalent to the fact that the restriction $d\a|_\xi$ to any
contact hyperplane is a non-degenerate 2-form. In particular, $\ker{}d\a$ is one-dimensional.

A vector field $X$ on $M$ is a \textit{contact vector field} if it
preserves the contact distribution $\xi$. 
In terms of contact forms this means that for every contact form $\a$, the Lie derivative
of $\a$ with respect to $X$ is proportional to $\a$:
\begin{equation}
\label{LieDerForm}
L_X\a=f_X\a
\end{equation}
where $f_X\in{}C^\infty(M)$.
The space of all contact vector fields (with complex coefficients) is a Lie algebra that we
denote $\CVect(M)$.

Let us now fix a contact form $\a$. A contact vector field $X$ is called 
strictly contact if it preserves $\a$, in other words, if $f_X=0$ everywhere on $M$. 
Strictly contact vector fields form a Lie subalgebra of $\CVect(M)$. 
There is one particular strictly contact vector field $Z$ called the \textit{Reeb field} (or
characteristic vector field).  It is
defined by the following two properties:
$$
Z\in\ker{}d\a,
\qquad
\a(Z)\equiv1.
$$

We will also consider the space, $\TVect(M)$, of (complex) vector fields \textit{tangent} to
the contact distribution.
This space is not a Lie subalgebra of $\Vect(M)$ that follows from non-integrability of
the contact distribution.

\subsection{The decomposition of $\Vect(M)$}

Let $\Vect(M)$ be the Lie algebra of all smooth vector fields (with complex coefficients)
on $M$.  The Lie bracket defines a natural action of $\CVect(M)$ on $\Vect(M)$.
In particular, the Lie bracket of a contact vector field with
a tangent vector field is again a tangent vector field.
Therefore, $\TVect(M)$ is a module over $\CVect(M)$.

\begin{prop}
\label{DecomPro}
The space $\Vect(M)$ is split into a direct sum of two $\CVect(M)$-modules:
$$
\Vect(M)\cong\CVect(M)\oplus\TVect(M).
$$
\end{prop}
\begin{proof}
Both spaces in the right hand side are $\CVect(M)$-modules. It then remains to check
that every vector field can be uniquely decomposed into a sum of a contact vector field and a tangent
vector field.

Given a vector field $X$, there exists a tangent vector field $Y$ such that $X-Y$ is contact. Indeed,
set $\b=L_X\a$ and consider the restriction of $\b|_\xi$ to a contact hyperplane. If $Y$ is a tangent
vector field then $L_Y\a=i_Y(d\a)$. Since $d\a$ in non-degenerate on $\xi$, then for any 1-form $\b$
there exists a tangent field $Y$ such that
$i_Y(d\a)|_\xi=\b|_\xi$. This means $X-Y$ is contact.

Furthermore, the intersection of $\CVect(M)$ and $\TVect(M)$ is zero. Indeed, let $X$ be a non-zero
vector field which is contact and tangent at the same time. Then $L_X\a=f\a$ for some function $f$ 
and $L_X\a=i_X(d\a)$. Since $\ker{}f\a$ contains $\xi=\ker\a$ while the restriction
$d\a|_\xi$ is non-degenerate, this is a contradiction.
\end{proof}

\section{The adjoint representation of $\CVect(M)$}

In this section we study the action of $\CVect(M)$ on itself.

\subsection{Fixing a contact form: contact Hamiltonians}\label{ConHamDepSec}

Let $M$ be orientable, fix a contact form $\a$ on $M$. Every contact vector field $X$ is then
characterized by a function:
$$
H=\a(X).
$$
This is a one-to-one correspondence between $\CVect(M)$ and the space $C^\infty(M)$ of
(complex valued) smooth functions on $M$, see, e.g.,
\cite{Arn}. We can denote the contact vector field corresponding to $H$ by $X_H$. The
function $H$ is called the contact Hamiltonian of $X_H$.

\begin{exe}
{\rm
The contact Hamiltonian of the Reeb field $Z$ is the constant function $H\equiv1$.
Note also that the function $f_X$ in (\ref{LieDerForm}) is given by the derivative
$f_{X_H}=Z(H)$.
}
\end{exe}

The Lie algebra $\CVect(M)$ is then
identified with
$C^\infty(M)$ equipped with the \textit{Lagrange bracket} defined by
$$
X_{\{H_1,H_2\}}:=
\left[
X_{H_1},X_{H_2}
\right].
$$
One checks that 
\begin{equation}
\label{adjrepForm}
\{H_1,H_2\}=X_{H_1}(H_2)-Z(H_1)\,H_2.
\end{equation}
The formula expresses the adjoint representation of $\CVect(M)$ in terms of contact Hamiltonians.
The second term in the right hand side shows that this action is different from the natural action of
$\CVect(M)$ on $C^\infty(M)$. Let us now clarify the geometric meaning of this action.

\subsection{Tensor densities on a contact manifold}\label{TenDenDefSec}

Let $M$ be an arbitrary smooth manifold of dimension $d$.
A \textit{tensor density} on $M$ of degree $\l\in\bbR$ is a section of
the line bundle $(\wedge^{d}T^*_\bbC{}M)^\l$. The space of $\l$-densities is denoted by
$\cF_\l(M)$.

Assume that $M$ is orientable and fix a volume form $\Omega$ on $M$.
This is a global section trivializing the above line bundle, so that $\cF_\l(M)$ can be
identified with $C^\infty(M)$. One then represents $\l$-densities in the form:
$$
\varphi=f\,\Omega^\l,
$$
where $f$ is a function.
\begin{exe}
{\rm
The space $\cF_0(M)\cong{}C^\infty(M)$ while the space $\cF_1(M)$ is nothing but the space of
differential $d$-forms. 
}
\end{exe}

If $M$ is compact then there is an invariant functional
\begin{equation}
\label{IntForm}
\int_M:\cF_1(M)\to\bbC.
\end{equation}
More generally, there is an invariant pairing
$$
\left\langle\cF_\l(M),\cF_{1-\l}(M)\right\rangle\to\bbC
$$
given by the integration of the product of tensor densities.

Let now $M$ be a contact manifold of dimension $d=2n+1$. In this case, there is another way
to define tensor densities. Consider the $(2n+2)$-dimensional submanifold
$S$ of the cotangent bundle $T^*M\setminus{}M$ that consists of all non-zero covectors 
vanishing on the contact distribution $\xi$. The restriction to $S$ of the canonical
symplectic structure on $T^*M$ defines a symplectic structure on $S$.
The manifold $S$ is called the \textit{symplectization} of $M$ (cf. \cite{Arn,AG}). Clearly 
$S$ is a line bundle over $M$, its sections are the 1-forms on $M$ vanishing 
on $\xi$.
Note that, in the case where $M$ is orientable, $S$ is a trivial line bundle
over $M$. 

There is a natural lift of $\CVect(M)$ to $S$. Indeed, a vector field $X$ on $M$ can
be lifted to $T^*M$, and, if $X$ is contact, then it preserves the subbundle
$S$. The space of sections $\mathrm{Sec}(S)$ is therefore a
$\CVect(M)$-module.

The sections of the bundle $S$ can be viewed as tensor 
densities of degree $\frac{1}{n+1}$ on $M$. 

\begin{prop}
\label{LemTDId}
There is a natural isomorphism of
$\CVect(M)$-modules
$$
\mathrm{Sec}(S)\cong\cF_{\frac{1}{n+1}}(M).
$$
\end{prop}
\begin{proof}
A section of $S$ is a 1-form on $M$ vanishing on the contact distribution.
For every contact vector field $X$ and a volume form $\Omega$ as in (\ref{CompNonint}) one
has
$$
L_X\Omega=(n+1)\,f_X\Omega.
$$
The Lie derivative of a tensor density of degree $\l$ is then given by
$$
L_X(f\,\Omega^\l)=
\left(
X(f)+\l(n+1)f_Xf
\right)
\Omega^\l.
$$
The result follows from formula (\ref{LieDerForm}).
\end{proof}

One can now represent $\l$-densities in terms of a contact form:
$
\varphi=f\,\a^{(n+1)\l}.
$

\subsection{Contact Hamiltonian as a tensor density}\label{CHSec}

In this section we identify the algebra $\CVect(M)$ with a space of tensor densities of
degree $-\frac{1}{n+1}$ on $M$; the adjoint action is nothing but a Lie derivative on this
space.  The result of this section is known (see \cite{Ovs} and \cite{OT}, Section 7.5) and
given here for the sake of completeness.

Let us define a different version of contact Hamiltonian of a contact vector field $X$ as a
$-\frac{1}{n+1}$-density on $M$:
$$
\cH:=\a(X)\,\a^{-1}.
$$
An important feature of this definition is that it is independent of the choice of $\a$.
Let us denote $X_\cH$ the corresponding contact vector field.

The space $\cF_{-\frac{1}{n+1}}(M)$ is now identified with $\CVect(M)$.
Moreover, the Lie bracket of contact vector fields corresponds to the Lie derivative.

\begin{prop}
\label{ConBraPro}
The adjoint representation of $\CVect(M)$ is isomorphic to
$\cF_{-\frac{1}{n+1}}(M)$.
\end{prop}
\begin{proof}
The Lagrange bracket coincides with a Lie derivative:
\begin{equation}
\label{LieDerPB}
\{\cH_1,\cH_2\}=L_{X_{\cH_1}}(\cH_2).
\end{equation}
This formula is equivalent to (\ref{adjrepForm}).
\end{proof}

Geometrically speaking, a contact Hamiltonian is not a function but rather a
tensor density of degree $-\frac{1}{n+1}$.

\subsection{Invariant functional on $\CVect(M)$}

Assume $M$ is compact and orientable, fix a contact form $\a$ and the corresponding volume
form $\Omega=\a\wedge{}d\a^n$.
The geometric interpretation of the adjoint action of $\CVect(M)$ implies the existence of an
invariant (non-linear) functional on $\CVect(M)$. 

Let $\CVect^*(M)$ be the set of contact vector fields
with nonvanishing contact Hamiltonians, this is an invariant open subset of $\CVect(M)$.

\begin{cor}
\label{InvFuncCorol}
The functional on $\CVect^*(M)$ defined by
$$
\cI:X_H\mapsto
\int_M{}H^{-(n+1)}\,\Omega
$$
is invariant. This functional is independent of the choice of the contact form.
\end{cor}
\begin{proof}
Consider is a contact vector
field $X_F$, then according to (\ref{adjrepForm}), one has
$$
L_{X_F}\,(H^{-(n+1)})
=
X_F\,(H^{-(n+1)})+(n+1)\,Z(F)
$$
so that the quantity $H^{-(n+1)}\,\Omega$ is a
well defined element of the space $\cF_1(M)$.
The functional $\cI$ is then given by the
invariant functional (\ref{IntForm}).

Furthermore, choose a different contact form $\a'=f\,\a$ and the corresponding volume form
$\Omega'=f^{n+1}\,\Omega$.
The contact Hamiltonian of the vector field $X_H$ with respect to the contact form $\a'$ is
the function
$H'=\a'(X_H)=f\,H$.
Hence, $H'^{-(n+1)}\,\Omega'=H^{-(n+1)}\,\Omega$ so that the functional $\cI$ is, indeed,
independent of the choice of the contact form.
\end{proof}

\section{The structure of $\TVect(M)$}

In this section we study the structure of
the space of tangent vector fields $\TVect(M)$ viewed as a $\CVect(M)$-module.

\subsection{A geometric realization}

Let us start with a geometric realization of the $\CVect(M)$-module structure on $\TVect(M)$
which if quite similar to that of Section \ref{CHSec}. 

Let $\Omega^2_0(M)$ be the space of 2-forms on $M$ vanishing on the
contact distribution. In other words, elements of $\Omega^2_0(M)$ are proportional to
$\a$:
$$
\omega=\a\wedge\b,
$$
where $\b$ is an arbitrary 1-form.

The following statement is similar to Proposition \ref{ConBraPro}.

\begin{thm}
\label{BizThm}
There is an isomorphism of $\CVect(M)$-modules
$$
\TVect(M)\cong\Omega^2_0(M)\otimes\cF_{-\frac{2}{n+1}}(M),
$$
where the tensor product is defined over $C^\infty(M)$.
\end{thm}
\begin{proof}
Let $M$ be orientable, fix a contact form $\a$ on $M$.
Consider a linear map
from $\TVect(M)$ to the space $\Omega^2_0(M)$ that associates to a tangent vector field $X$
the 2-form
$$
\left\langle
X,\a\wedge{}d\a
\right\rangle=
-\a\wedge{}i_Xd\a.
$$
This map is bijective since the restriction $d\a|_\xi$ of the 2-form $d\a$ to the
contact hyperplane $\xi$ is non-degenerate. 

However, the above map depends on the choice of the contact form and, therefore, cannot be
$\CVect(M)$-invariant.
In order to make this map independent of the choice of $\a$, one defines the following map
\begin{equation}
\label{GeomRealForm}
X\mapsto
\left\langle
X,\a\wedge{}d\a
\right\rangle\otimes\a^{-2}
\end{equation}
with values in $\Omega^2_0(M)\otimes\cF_{-\frac{2}{n+1}}(M)$.
Note that the term $\a^{-2}$ in the right-hand-side is a well defined element of the
space of tensor densities
$\cF_{-\frac{2}{n+1}}(M)$, see Section \ref{TenDenDefSec}.

It remains to check the $\CVect(M)$-invariance of the map (\ref{GeomRealForm}).
Let $X_H$ be a contact vector field, one has
$$
\begin{array}{rcl}
L_{X_H}\left(
\left\langle
X,\a\wedge{}d\a
\right\rangle\otimes\a^{-2}
\right)&=&
\left\langle
[X_H,X],\a\wedge{}d\a
\right\rangle\otimes\a^{-2}\\[6pt]
&&+
\left\langle
X,f_X\a\wedge{}d\a+\a\wedge{}df_X\a
\right\rangle\otimes\a^{-2}
-
\left\langle
X,\a\wedge{}d\a
\right\rangle\otimes{}(2f_X\a^{-2})\\[6pt]
&=&
\left\langle
[X_H,X],\a\wedge{}d\a
\right\rangle\otimes\a^{-2}.
\end{array}
$$
Hence the result.
\end{proof}

The isomorphism (\ref{GeomRealForm}) identifies the $\CVect(M)$-action on $\TVect(M)$ by
Lie bracket with the usual Lie derivative. It is natural to say that this map
defines an analog of contact Hamiltonian of a tangent vector field.

\subsection{A skew-symmetric pairing on $\TVect(M)$ over $\CVect(M)$}

There exists an invariant skew-symmetric bilinear map
from $\TVect(M)$ to $\CVect(M)$ that can be understood as a ``symplectic structure'' on the
space $\TVect(M)$ over $\CVect(M)$.

\begin{thm}
\label{BilMapThm}
There exists a non-degenerate skew-symmetric invariant bilinear map
$$
\cB:\TVect(M)\wedge\TVect(M)\to\CVect(M),
$$
where the $\wedge$-product is defined over $C^\infty(M)$.
\end{thm}
\begin{proof}
Assume first that $M$ is orientable and fix the contact form $\a$. 
Given $2$ tangent vector fields $X$ and $Y$, consider the function
$$
H_{X,Y}=
\left\langle
X\wedge{}Y\,,\,
d\a
\right\rangle.
$$
Define first a $(2n)$-linear map $B$ from $\TVect(M)$ to $C^\infty(M)$ by
\begin{equation}
\label{BilMapDepdForm}
B_\a:X\wedge{}Y\mapsto{}H_{X,Y}.
\end{equation}
The definition of the function $H_{X,Y}$ and thus of the map $B_\a$ depends on
the choice of $\a$. Our task is to understand it as a map with values in
$\CVect(M)$ which is independent of the choice of the contact form. 
This will, in particular, extend the definition to the case where
$M$ is not orientable.

It turns out that the above function $H_{X,Y}$ is a well-defined contact Hamiltonian.

\begin{lem}
\label{ConstLem}
Choose a different contact form $\a'=f\,\a$, then $H'_{X,Y}=f\,H_{X,Y}$.
\end{lem}
\begin{proof}
By definition,
$$
H'_{X,Y}=
\left\langle
X\wedge{}Y,d\a'
\right\rangle=
f\left\langle
X\wedge{}Y,d\a
\right\rangle
+\left\langle
X\wedge{}Y,df\wedge\a
\right\rangle=
f\,H_{X,Y}
$$
since the second term vanishes.
\end{proof}

We observe that the function $H_{X,Y}$ depends on
the choice of $\a$ precisely in the same way as a contact Hamiltonian 
(cf. Section \ref{ConHamDepSec}). 
It follows that the bilinear map
\begin{equation}
\label{BilMapForm}
\cB:X\wedge{}Y\mapsto{}H_{X,Y}\,\a^{-1}
\end{equation}
with values in $\cF_{-\frac{1}{n+1}}\cong\CVect(M)$ (cf. Section \ref{CHSec})
is well-defined and independent of the choice of $\a$.

It remains to check that the constructed map (\ref{BilMapForm}) is $\CVect(M)$-invariant.
This can be done directly but also follows from

\begin{prop}
\label{BilMapAltPro}
The Lie bracket of two tangent vector fields $X,Y\in\TVect(M)$ is of the form
\begin{equation}
\label{BilMapCommForm}
[X,Y]=
\cB(X,Y)
\quad
+
\quad
\hbox{(tangent vector field)}
\end{equation}
\end{prop}
\noindent
\begin{proof}
Consider the decomposition from Proposition \ref{DecomPro} applied to the Lie bracket $[X,Y]$.
The ``non-tangent'' component of $[X,Y]$
is a contact vector field with contact Hamiltonian $\a([X,Y])$. One has
$$
i_{[X,Y]}\a=
\left(L_X\,i_Y-i_Y\,L_X\right)\a=
-i_Y\,L_X\,\a=
-i_Y\,i_X\,d\a=
H_{X,Y}
$$
The result follows.
\end{proof}
Theorem \ref{BilMapThm} is proved.
\end{proof}
Proposition \ref{BilMapAltPro} is an alternative definition of $\cB$:
the map $\cB$ measures the failure of the Lie bracket of two
tangent vector fields to remain tangent.

\section{Heisenberg structures}

In order to investigate the structure of $\TVect(M)$ as a $\CVect(M)$-module in more
details, we will write explicit formul{\ae} for the $\CVect(M)$-action.

We assume that there is an action of the Heisenberg Lie algebra $\fh_n$ on $M$,
such that the center acts by the Reeb field while the generators are tangent to the contact
structure. We then say that $M$ is equipped with the Heisenberg structure.
Existence of a globally defined Heisenberg structure is a strong condition on $M$, however,
locally such structure always exists.

\subsection{Definition of a Heisenberg structure}

Recall that the Heisenberg Lie algebra $\fh_n$ is a nilpotent Lie algebra of dimension
$2n+1$ with the basis
$
\left\{a_1,\ldots,a_{n},b_1,\ldots,b_{n},z\right\}
$
and the commutation relations
$$
\left[
a_i,b_j
\right]=\d_{ij}\,z,
\qquad
\left[
a_i,a_j
\right]=
\left[
b_i,b_j
\right]=
\left[
a_i,z
\right]=
\left[
b_i,z
\right]=0,
\qquad
i,j=1,\ldots,n.
$$
The element $z$ spans the one-dimensional center of $\fh_n$.

\begin{rmk}
{\rm
The algebra $\fh_n$ naturally appears in the context of symplectic geometry as a Poisson
algebra of linear functions on the standard $2n$-dimensional symplectic space.
}
\end{rmk}

We say that $M$ is equipped with a {\it Heisenberg structure} if one fixes a contact form $\a$
on $M$ and a $\fh_n$-action spanned by $2n+1$ vector fields
$
\left\{A_1,\ldots,A_{n},B_1,\ldots,B_{n},Z\right\},
$
such that the $2n$ vector fields $A_i,B_j$ are independent at any point and tangent to the
contact structure:
$$
i_{A_i}\,\a=i_{B_j}\,\a=0
$$
and $[A_i,B_i]=Z$, where $Z$ is the Reeb field, while the other Lie brackets are zero.

\subsection{Example: the local Heisenberg structure}\label{DarSect}

The Darboux theorem states that locally contact manifolds are diffeomorphic to each other.
An effective way to formulate this theorem is to say that in a neighborhood of any point of
$M$ there is a system of local coordinates $(x_1,\ldots,x_{n},y_1,\ldots,y_{n},z)$ such that
the contact structure $\xi$ is given by the 1-form
$$
\a=\sum_{i=1}^{n}
\frac{x_i\,dy_i-y_i\,dx_i}{2}+dz.
$$
These coordinates are called the Darboux coordinates.

\begin{prop}
\label{HeisProp}
The vector fields
\begin{equation}
\label{DbarFields}
A_i=
\frac{\partial}{\partial x_i}+
\frac{y_i}{2}\,\frac{\partial}{\partial z},
\qquad
B_i=
-\frac{\partial}{\partial y_i}+
\frac{x_i}{2}\,\frac{\partial}{\partial z},
\qquad
Z=
\frac{\partial}{\partial z},
\end{equation}
where $i=1,\ldots,n$, define a Heisenberg structure on $\bbR^{2n+1}$.
\end{prop}

\begin{proof}
One readily checks that $A_i,B_j$ are tangent and
$$
\left[
A_i,B_j
\right]=\d_{ij}\,Z
$$
while other commutation relations are zero. The vector field $Z$ is nothing but the Reeb
field.
\end{proof}

There is a well-known formula for a contact vector field in the Darboux
coordinates (see, e.g., \cite{Arn,AG,Kir}). 
We will not use this formula since the expression in terms of the Heisenberg structure is
much simpler.

\subsection{Contact vector fields and Heisenberg structure}

Assume that $M$ is equipped with an arbitrary Heisenberg structure. It turns out that every
contact vector fields can be expressed in terms of the basis of the $\fh_n$-action by a
universal formula. 

\begin{prop}
Given an arbitrary Heisenberg structure on $M$, a contact vector field with a contact
Hamiltonian $H$ is given by the formula
\begin{equation}
\label{DecomConForm}
X_H=
H\,Z-
\sum_{i=1}^n\left(
A_i(H)\,B_i-B_i(H)\,A_i
\right).
\end{equation}
\end{prop}
\begin{proof}
Let us first check that the vector field (\ref{DecomConForm}) is, indeed, contact.
If $X$ be as the right-hand-side of (\ref{DecomConForm}), then the Lie derivative
$L_X\a:=\left(d\circ{}i_X+i_X\circ{}d\right)\a$ is given by
$$
L_X\a=
dH-
\sum_{i=1}^n\left(
A_i(H)\,i_{B_i}-B_i(H)\,i_{A_i}
\right)d\a.
$$
To show that the 1-form $L_X\a$ is proportional to $\a$, it suffice to check that
$$
i_{A_i}\left(L_X\a\right)=i_{B_j}\left(L_X\a\right)=0
\qquad
\hbox{for all}
\quad
i,j=1,\ldots{}n.
$$
The first relation is a consequence of the formul{\ae} $i_{A_i}\left(dH\right)=A_i(H)$
together with
\begin{equation}
\label{PdstavForm}
i_{A_i}i_{B_j}\,d\a
=i_{A_i}\left(L_{B_j}\a\right)
=i_{[A_i,B_j]}\a
=\d_{ij}\,i_{Z}\,\a
=\d_{ij},
\qquad
i_{A_i}i_{A_j}\,d\a=i_{B_i}i_{B_j}\,d\a=0.
\end{equation}
The second one follows from the similar relations for
$i_{B_j}$.

Second, observe that, if $X$ be as in (\ref{DecomConForm}), then $i_X\,\a=H$. 
This means that the
contact Hamiltonian of the contact vector field (\ref{DecomConForm}) is precisely $H$.
\end{proof}

Note that a formula similar to (\ref{DecomConForm}) was used in \cite{Kir} to define a
contact structure.

\subsection{The action of $\CVect(M)$ on $\TVect(M)$}

Since $2n$ vector fields $A_i$ and $B_j$ are linearly independent at
any point, they form a basis of $\TVect(M)$ over $C^\infty(M)$. 
Therefore, an arbitrary tangent vector field $X$ has a unique decomposition
\begin{equation}
\label{DecomArbForm}
X=
\sum_{i=1}^n\left(
F_i\,A_i+ G_i\,B_i
\right),
\end{equation}
where $(F_i,G_j)$ in an $2n$-tuple of smooth functions on $M$.
The space $\TVect(M)$ is now identified with the direct sum
$$
\TVect(M)\cong
\underbrace{C^\infty(M)\oplus\cdots\oplus{}C^\infty(M)}_{\mbox{{\tiny $2n$ times}}},
$$
Let us calculate explicitly the action of $\CVect(M)$ on $\TVect(M)$.

\begin{prop}
\label{ActionExpProp}
The action of $\CVect(M)$ on $\TVect(M)$ is given by the first-order
$(2n\times{}2n)$-matrix differential operator
\begin{equation}
\label{MatGamForm}
X_H
\left(
\begin{array}{c}
F\\[8pt]
G
\end{array}
\right)
=
\left(
X_H\cdot{\mathbf{1}}
-
\left(
\begin{array}{rr}
AB(H) & BB(H)\\[8pt]
-AA(H) & -BA(H)
\end{array}
\right)
\right)
\left(
\begin{array}{c}
F\\[8pt]
G
\end{array}
\right)
\end{equation}
where $F$ and $G$ are $n$-vector functions, $\mathbf{1}$ is the unit $(2n\times{}2n)$-matrix,
$AA(H),AB(H),BA(H)$ and
$BB(H)$ are $(n\times{}n)$-matrices, namely
$$
AA(H)_{ij}=A_iA_j(H),
$$
the three other expressions are similar.
\end{prop}

\begin{proof}
Straightforward from (\ref{DecomConForm}) and (\ref{DecomArbForm}).
\end{proof}

\begin{prop}
\label{BilExpProp}
The bilinear map (\ref{BilMapDepdForm}) has the following explicit expression:
$$
H_{X,\widetilde{X}}=
\sum_{i=1}^n
\left|
\begin{array}{cc}
F_i & \widetilde{F}_i\\
G_i & \widetilde{G}_i
\end{array}
\right|,
$$
where 
$X=\sum_{i=1}^n
(F_i\,A_i+ G_i\,B_i),$ 
and 
$\widetilde{X}=\sum_{j=1}^n
(\widetilde{F}_j\,A_j+ \widetilde{G}_j\,B_j)$.
\end{prop}
\begin{proof}
This follows from the definition (\ref{BilMapDepdForm}) and formula (\ref{PdstavForm}).
\end{proof}
\noindent
Note that formula (\ref{MatGamForm}) implies that $H_{X,\widetilde{X}}$ transforms as a
contact Hamiltonian according to (\ref{adjrepForm}) since the partial traces of the
$(2n\times{}2n)$-matrix in (\ref{MatGamForm}) are $A_iB_i(H)-B_iA_i(H)=Z(H)$.

\vskip 0.5cm

\textbf{Acknowledgments}.
I am grateful to C. Duval and S. Tabachnikov for their interest in this work and a careful
reading of a preliminary version of this paper.

\vskip 0.5cm


\end{document}